\documentclass[12pt]{amsart}

\newtheorem{theorem}{Theorem}[section]
\newtheorem{lemma}[theorem]{Lemma}

\newtheorem{conjecture}[theorem]{Conjecture}

\theoremstyle{plain}
\newtheorem{definition}[theorem]{Definition}
\newtheorem{example}[theorem]{Example}

\newtheorem{question}[theorem]{Question}

\newtheorem{setting}[theorem]{Setting}

\theoremstyle{definition}

\theoremstyle{remark}

\numberwithin{equation}{section}

\usepackage{fancyhdr}
\usepackage{amscd}
\usepackage{amsmath}
\usepackage{amsthm}
\usepackage{amssymb}

\usepackage[all]{xy}

\begin{document}

\title[]{Deformations of Chow groups via Cyclic homology}

 \author{Sen Yang \\ \\ Dedicated to Jinke Hai with profound gratitude }

 \address{School of Mathematics and Finance \\  Chuzhou University \\
 Chuzhou, 239000, China\\
 }

\address{Applied Mathematics Research Center \\ Chuzhou University \\
 Chuzhou, 239000, China \\
}

 \email{yangsen.2014@tsinghua.org.cn, 101012424@seu.edu.cn}

\date{}

\maketitle

\begin{abstract} 
Let $X$ be a smooth projective variety over an arbitrary field $k$ of characteristic zero. We explore infinitesimal deformations of the Chow group $CH^{p}(X)$ via its formal completion $\widehat{CH}^{p}$, a functor defined on the category of local augmented Artinian $k$-algebras. Under a natural vanishing condition on Hodge cohomology groups, for certain augmented graded Artinian $k$-algebras $A$ with the maximal ideal $m_{A}$, we prove that
\[
\widehat{CH}^{p}(A) \cong H^{p}(X, \Omega^{p-1}_{X/ k})\otimes_{k}m_{A}.
\]This extends earlier results of Bloch and others from the case where $k$ is algebraic over $\mathbb{Q}$ to arbitrary fields of characteristic zero, and gives a partial affirmative answer to a general question linking the pro-representability of Chow groups to a specific set of Hodge-theoretic vanishing conditions.  


\end{abstract}

\section{Introduction}

Let $X$ be a smooth projective variety over an arbitrary field $k$ of characteristic zero. For an integer $p$ with $1 \leq p \leq \mathrm{dim}(X)$, let $CH^{p}(X)$ denote the Chow group of codimension-$p$ algebraic cycles modulo rational equivalence.  As a fundamental invariant encoding the geometric and arithmetic information of algebraic varieties, the Chow group plays a central role in algebraic geometry, linking to various topics such as intersection theory and Hodge theory. This paper is devoted to the study of infinitesimal deformations of the Chow group $CH^{p}(X)$.

Let $Art_{k}$ denote the category of local Artinian $k$-algebras $A$ with the maximal ideal $m_{A}$ such that $A/m_{A}=k$. For $A \in  Art_{k}$, let $X_A=X\times_{\mathrm{Spec}(k)}\mathrm{Spec}(A)$ be the infinitesimal deformation of $X$ over $A$. We denote by $K^{M}_{p}(O_{X_{A}})$ the Milnor K-theory sheaf associated to the following presheaf
\[
 U \to K^{M}_{p}(O_{X_{A}}(U)),
\]where $U=\mathrm{Spec}(R) \subset X$ is open affine, $K^{M}_{p}(O_{X_{A}}(U))=K^{M}_{p}(R \otimes_{k}A)$. When $A=k$, we write $K^{M}_{p}(O_{X})=K^{M}_{p}(O_{X_{k}})$.

For an abelian group $M$, let $M_{\mathbb{Q}}$ denote $M \otimes_{\mathbb{Z}} \mathbb{Q}$ in the sequel. We define the formal completion of $CH^{p}(X)$ as the following functor (cf. Bloch \cite{Bl3})
{\footnotesize
\begin{align} \label{eq:formalcom}
\widehat{CH}^{p}: \ A \longmapsto  \mathrm{ker} \ \{ H^{p}(X, K^{M}_{p}(O_{X_A}))_{\mathbb{Q}}  \xrightarrow{aug} H^{p}(X, K^{M}_{p}(O_{X}))_{\mathbb{Q}} \},
\end{align}
}where $A \in  Art_{k}$, and the map $aug$ is induced by the augmentation $A \to k$.

The Bloch formula (cf. \cite{Bl2-Annals, Ka,Kerz,Quillen,Soule}) states
\begin{equation}
CH^{p}(X)_{\mathbb{Q}}\cong H^{p}(X, K^{M}_{p}(O_{X}))_{\mathbb{Q}},
\end{equation}
meaning that the cohomology group $H^{p}(X, K^{M}_{p}(O_{X}))_{\mathbb{Q}}$ identifies with the Chow group. Therefore, $\widehat{CH}^{p}(A)$, as defined above, can capture useful information about infinitesimal deformations of $CH^{p}(X)_{\mathbb{Q}}$.

Pro-representable functors, recalled in Definition \ref{d:rep}, are powerful tools for studying deformation problems. We are interested in the pro-representability of the functor $\widehat{CH}^{p}$. The case of $p=1$ is classical, cf. Lemma \ref{l:ch1} below. When $p \geq 2$, the study of the pro-representability of the functor $\widehat{CH}^{p}$ remains a challenge. Pioneering work by Bloch, Stienstra, Mackall and others has shown that the pro-representability of Chow groups is influenced by the Hodge theory of $X$. When the ground field $k$ is algebraic over the rational number field $\mathbb{Q}$, we have the following result.
\begin{theorem} [\cite{Bl3, Ma,Y6}]\label{t:b-m-yang}
Let $X$ be a smooth projective variety over a field $k$ of dimension $d \geq 2$, where $k$ is an algebraic field extension of $\mathbb{Q}$. Let $p$ be an integer such that $2 \leq p \leq d$. We suppose that $X$ satisfies the following condition:
\begin{equation} \label{eq: condi-vanish1}
H^{p}(X, \Omega^{i}_{X/k})= H^{p+1}(X, \Omega^{i}_{X/k})= \cdots =H^{2p-1-i}(X,\Omega^{i}_{X/k})=0,
\end{equation}where $i$ is an integer such that $0 \leq i \leq p-2$. Then, for any $A \in Art_{k}$, there is an isomorphism
\begin{equation} \label{eq:mainiso}
\widehat{CH}^{p}(A)\cong H^{p}(X, \Omega^{p-1}_{X/ k})\otimes_{k}m_{A},
\end{equation}where $m_{A}$ is the maximal ideal of $A$. Consequently, the functor $\widehat{CH}^{p}$ is pro-representable.

\end{theorem}

The cases $p=2$ and $p=3$ were proved by Bloch \cite{Bl3} and by Mackall \cite{Ma} respectively. We generalize their results into the above format, cf. Theorem 4.10 in \cite{Y6}. The study of $\widehat{CH}^{2}$ motivated Bloch to propose the following important conjecture, which remains one of the major questions concerning Chow groups and has been intensively investigated (cf. \cite{BKL,BS,Hu,PW,V2,V3}).
\begin{conjecture} [\cite{Bl3}] \label{c:Bloch}
Let $X$ be a connected smooth complex projective surface with trivial geometric genus.\footnote{Recall that the geometric genus $p_{g}(X)=\mathrm{dim} H^{0}(X,\Omega^{2}_{X/k})=\mathrm{dim} H^{2}(X,O_{X})$.} Then the Albanese map
\[
CH^2_{\text{deg } 0}(X) \to Alb(X),
\]
is an isomorphism. Here $CH^2_{\text{deg } 0}(X)$ is the subgroup of $CH^{2}(X)$ consisting of zero cycles with degree zero and $Alb(X)$ is the Albanese variety.
\end{conjecture}

If the ground field $k$ contains transcendental elements over $\mathbb{Q}$, the non-trivial absolute differentials $\Omega^{1}_{k/ \mathbb{Q}}$ may appear. This results in the difficulty in computing $\widehat{CH}^{p}(A)$. When $p=2$, this subtlety is handled by Stienstra.

\begin{theorem}[\cite{Stien1}] \label{t:stien}
Let X be a smooth projective surface over a field $k$ which is an extension of $\mathbb{Q}$ of finite transcendence degree. Then the
following are equivalent:

$\mathrm{(1)}$ The map $H^{1}(X,\Omega^{1}_{X/k}) \to H^{2}(X,O_{X})\otimes_{k} \Omega^{1}_{k/ \mathbb{Q}}$ which is induced by Gauss-Manin connection is surjective.

 $\mathrm{(2)}$ For $A \in Art_{k}$, the functor assigning $A$ to the kernel of the map
\[
 \widehat{CH}^{2}(A) \to H^{2}(\Omega^{1}_{X/k})\otimes_{k} m_{A},
\]is naturally isomorphic to the one which assigns $A$ to the group $H^{2}(O_{X})\otimes_{k} \frac{\Omega^{1}_{A/k}}{dm_{A}}$.

\end{theorem}

When the ground field is the complex number field $\mathbb{C}$, the Bloch-Beilinson conjecture predicts that there is a conjectural filtration of the Chow group $CH^{p}(X)_{\mathbb{Q}}$
\begin{align*}
CH^{p}(X)_{\mathbb{Q}} & =F^{0}CH^{p}(X)_{\mathbb{Q}} \supset F^{1}CH^{p}(X)_{\mathbb{Q}} \\
& \supset \cdots \supset F^{p}CH^{p}(X)_{\mathbb{Q}} \supset F^{p+1}CH^{p}(X)_{\mathbb{Q}}=0,
\end{align*}whose graded quotients
\[
Gr^{i}CH^{p}(X)_{\mathbb{Q}} =\dfrac{F^{i}CH^{p}(X)_{\mathbb{Q}}}{F^{i+1}CH^{p}(X)_{\mathbb{Q}}}, \ \mathrm{where} \ i =0, \ \cdots, \ p,
\]admits a Hodge-theoretic description.

An affirmative answer to this conjecture would greatly help us understand the structure of the Chow group $CH^{p}(X)_{\mathbb{Q}}$ and have a far-reaching influence, cf. \cite{Ja, Le, Mu, Sa}.  The first two steps are known and $F^{2}CH^{p}(X)_{\mathbb{Q}}$ is the kernel of the Deligne cycle class map
\[
r: CH^{p}(X)_{\mathbb{Q}} \to  H_{\mathcal{D}}^{2p}(X, \mathbb{Z}(p)).
\]We refer the readers to \cite{El-Z, EV, Jannsen} for the definition of the Deligne cohomology $H_{\mathcal{D}}^{2p}(X, \mathbb{Z}(p))$  and the construction  of the Deligne cycle class.

In \cite{GG}, Green and Griffiths propose a definition of such a filtration in a Hodge-theoretic way. Of particular relevance to our study of Chow groups is the relationship between the transcendence degree of the ground field $k$ and the Chow groups themselves.

\begin{conjecture} [cf. Implication 1.2 on page 478 of \cite{GG}] \label{c:BB-number}
If $X$ is a smooth projective variety over a number field $k$, then
\[
F^{2}CH^{p}(X)_{\mathbb{Q}}=0,
\]
where $F^{2}CH^{p}(X)_{\mathbb{Q}}$ is the filtration induced from $F^{2}CH^{p}(X_{\mathbb{C}})_{\mathbb{Q}}$ under the natural map
\[
CH^{p}(X) \to CH^{p}(X_{\mathbb{C}}),
\] where $X_{\mathbb{C}}:= X \times_{\mathrm{Spec}(k)} \mathrm{Spec}(\mathbb{C})$. 

\end{conjecture}

For the ground field of transcendence degree one, extending earlier examples by Bloch, Nori and Schoen, Green-Griffiths-Paranjape \cite{GGP} has proved that there exists a variety $V$ such that $F^{2}CH^{p}(V)_{\mathbb{Q}}$ is nontrivial.

This demonstrates an interesting idea: the transcendence degree of the ground field indeed influences the structure of Chow groups. To combine this idea with the work by Bloch, Stienstra, Mackall and the author on the pro-representability of Chow groups, we investigate the influence of the transcendence degree of the ground field on the pro-representability of the functor $\widehat{CH}^{p}$, defined in \eqref{eq:formalcom}, and focus on the following question in this paper.

\begin{question}  \label{q:main}
Let $X$ be a smooth projective variety of dimension $d \geq 2$ over an arbitrary field $k$ of characteristic zero. Let $p$ be an integer such that $2 \leq p \leq d$. If $X$ satisfies the condition \eqref{eq: condi-vanish1} in Theorem \ref{t:b-m-yang}, does the isomorphism \eqref{eq:mainiso} still hold?
\end{question}

The principal result of this paper is to partially answer this question, by partially generalizing Theorem \ref{t:b-m-yang} to the case when the ground field $k$ is an arbitrary field.
\begin{theorem}[cf. Theorem \ref{t:yang1} below]  \label{t:yang2}
Let $X$ be a smooth projective variety over an arbitrary field $k$ of characteristic zero. If $X$ satisfies the condition \eqref{eq: condi-vanish1}, then for $A \in Art_{k}$ which is also a graded $k$-algebra $A=\oplus_{l \geq 0}A_{l}$ such that $A_{0}=k$, there is an isomorphism
\[
\widehat{CH}^{p}(A) \cong H^{p}(X, \Omega^{p-1}_{X/ k})\otimes_{k}m_{A}.
\]

\end{theorem}

Our proof leverages the powerful machinery of cyclic homology, following the philosophy of Goodwillie \cite{GoodW, Good} that relates relative K-theory to relative cyclic homology. We establish a series of technical lemmas that allow us to control the relevant cohomology groups of sheaves of differentials, ultimately reducing the computation of $\widehat{CH}^{p}(A)$ to a single Hodge cohomology group of $X$.

This paper is structured as follows. In Section 2, we recall the necessary background on Hochschild and cyclic homology. Section 3 is devoted to the proof of our main result, Theorem \ref{t:yang1}.

\textbf{Notations:} 

(1) We denote by $k$ an arbitrary  field of characteristic zero, and denote by $Art_{k}$ the category of local Artinian $k$-algebras $A$ with the maximal ideal $m_{A}$ such that $A/m_{A}=k$.

(2) For any abelian group $M$, $M_{\mathbb{Q}}$ denotes $M \otimes_{\mathbb{Z}} \mathbb{Q}$. 

\section{Cyclic homology}

Let $R$ be an associative commutative unital $k$-algebra. The Hochschild complex (with coefficients in $R$) is 
\begin{equation*}
C(R): \ \ \cdots \to R^{\otimes (n+1)} \xrightarrow{b} R^{\otimes n} \to \cdots \to R \otimes_{k}R \xrightarrow{b} R,
\end{equation*} in which $n$ is a non-negative integer and $R^{\otimes (n+1)}$ is in  degree $n$. The map $b: \ R^{\otimes (n+1)} \xrightarrow{b} R^{\otimes n}$ is the Hochschild boundary defined as (cf. section 1.1.1 of \cite{Loday})
\begin{align*}
b: a_{0} \otimes a_{1} \otimes  a_{2} \otimes \cdots \otimes a_{n}  \to & \sum_{i=1}^{n-1}(-1)^{i} a_{0} \otimes a_{1} \otimes  \cdots \otimes a_{i} a_{i+1} \otimes  \cdots \otimes a_{n} \\
& + a_{0}a_{1} \otimes a_{2} \otimes  \cdots \otimes a_{n} \\
& +(-1)^{n} a_{n}a_{0} \otimes a_{1} \otimes  \cdots \otimes a_{n-1}.
\end{align*}

The Hochschild homology group $HH_{n}(R)$ is defined to be the homology of the Hochschild complex $C(R)$
\begin{equation}  \label{eq:hhalg}
HH_{n}(R):=H_{n}(C(R)).
\end{equation}

To illustrate this definition, we compute $HH_{0}(R)$ and $HH_{1}(R)$ by considering the following two boundary maps
\[
 \cdots \to R \otimes_{k}R \otimes_{k}R \xrightarrow{b} R \otimes_{k}R \xrightarrow{b} R.
\]

The map $b: R \otimes_{k}R \to R$ is 
\begin{equation}  \label{eq:map b1}
a_{0}\otimes a_{1} \to a_{0}a_{1}-a_{1}a_{0}, 
\end{equation}whose image generates the commutator $[R,R]$ of $R$. Hence, we have 
\[
HH_{0}(R)=R/[R,R].
\]Since $R$ is commutative here, $[R,R]$ is trivial, we have
\[
HH_{0}(R)=R.
\]
 
The map $b: R \otimes_{k}R \otimes_{k}R \to R \otimes_{k}R$ is 
\begin{equation}  \label{eq:map b2}
a_{0}\otimes a_{1}\otimes a_{2} \to a_{0} a_{1}\otimes a_{2}-a_{0}\otimes a_{1} a_{2} + a_{2}a_{0}\otimes a_{1}
\end{equation}

Since the map \eqref{eq:map b1} is trivial, $HH_{1}(R)$ is the quotient of $R \otimes_{k}R$ modulo the imgae of the map \eqref{eq:map b2}. Consequently, there is a well-defined map
\begin{align*}
HH_{1}(R) \to \Omega^{1}_{R/k} \\
a\otimes b \to adb,
\end{align*}which is indeed an isomorphism. 

\begin{lemma} [cf. Prop. 1.1.10 of \cite{Loday}]  \label{l:HH1}
With the notation as above, there is a canonical isomorphism
\[
HH_{1}(R) \cong \Omega^{1}_{R/k}.
\]

\end{lemma}

This uncovers the connection between the first Hochschild homology and differential forms. The Hochschild-Kostant-Rosenberg (HKR) theorem asserts that, for smooth algebras, Hochschild homology coincides with 
differential forms.
\begin{theorem} [HKR theorem]  \label{l:HKR-smooth}
If $R$ is smooth over $k$, then there is a canonical isomorphism
\[
HH_{n}(R) \cong \Omega^{n}_{R/k}.
\]

\end{theorem}

The cyclic bicomplex $CC(R)$ is a first-quadrant bicomplex
\[
  \begin{CD}
      \vdots @<<<  \vdots @<<<  \vdots @<<< \vdots  @<<< \\ 
       @VbVV @V-b'VV @VbVV @V-b'VV @.\\
      R^{\otimes 3}  @<1-t<< R^{\otimes 3} @<N<<  R^{\otimes 3}  @<1-t<<  R^{\otimes 3}  @<N<<  \\
     @VbVV @V-b'VV @VbVV @V-b'VV @.\\
      R^{\otimes 2}  @<1-t<< R^{\otimes 2} @<N<<  R^{\otimes 2}  @<1-t<<  R^{\otimes 2}  @<N<<  \\
    @VbVV @V-b'VV @VbVV @V-b'VV @.\\
      R  @<1-t<<  R @<N<<  R  @<1-t<<  R  @<N<<,
  \end{CD}
\] in which the module $R$ in the left corner is of bidegree $(0,0)$. Here the horizontal differentials $t$ and $N$ are alternatively cyclic operators and norm operators, and the vertical complexes are either Hochschild complex $C(R)$ or its variant, cf. section 5.1.1 of \cite{Loday}.

The cyclic homology of $R$ is defined to be the homology group of the total complex $\mathrm{Tot}(CC(R))$ of the bicomplex $CC(R)$
\begin{equation}  \label{eq:hcalg}
HC_{n}(R):=H_{n}(\mathrm{Tot}(CC(R))).
\end{equation} We refer to section 5.1 of \cite{Loday} for details. Classically, we have
\begin{equation*}
HC_{0}(R)=HH_{0}(R)=R.
\end{equation*}

\begin{example} [cf. Ex 4.1.8 of \cite{Loday}] \label{ex:HC-k}
The cyclic homology of $k$ can be computed as 
\begin{align*}
\begin{cases}
 \begin{CD}
   & HC_{n}(k)=k,  n \ \mathrm{is \ even}, \\
  &HC_{n}(k)=0, n \ \mathrm{is \ odd}.
 \end{CD}
\end{cases}
\end{align*}
\end{example}

The $n$-th de Rham cohomology of $R$, denoted $H^{n}_{dR}(R)$, is defined to be the  $n$-th cohomology of the following de Rham complex
\begin{equation*} \label{eq: derham}
\cdots \to \Omega^{n-1}_{R/k} \xrightarrow{d} \Omega^{n}_{R/k} \xrightarrow{d} \Omega^{n+1}_{R/k} \to \cdots.
\end{equation*}
More precisely, let $\mathrm{Ker}(d)$ denote the kernel of $\Omega^{n}_{R/k} \xrightarrow{d} \Omega^{n+1}_{R/k}$, and let $d\Omega^{n-1}_{R/k}$ denote the image of $\Omega^{n-1}_{R/k} \xrightarrow{d} \Omega^{n}_{R/k}$. Then the de Rham cohomology $H^{n}_{dR}(R)$ is defined as a quotient
\begin{equation*} \label{eq:derhamgp}
H^{n}_{dR}(R):=\frac{\mathrm{Ker}(d)}{d\Omega^{n-1}_{R/k}}.
\end{equation*}

Loday and Quillen proved that de Rham cohomologies of smooth algebras are related to cyclic homology groups. 
\begin{lemma} [\cite{LQ}] \label{l:cyc-hom-smooth}
If $R$ is smooth over $k$, then there is a canonical isomorphism
\[
HC_{n}(R) \cong \frac{\Omega^{n}_{R/k}}{d\Omega^{n-1}_{R/k}} \oplus H^{n-2}_{dR}(R) \oplus H^{n-4}_{dR}(R)\oplus \cdots.
\]The last summand is $H^{0}_{dR}(R)$ or $H^{1}_{dR}(R)$, depending on whether $n$ is even or odd respectively.

\end{lemma}

\begin{theorem}[cf. Th 2.2.1 of \cite{Loday}] \label{t:SBI}
For an associative commutative unital $k$-algebra $R$, there is a natural long exact 
sequence
\begin{align} \label{eq:SBI}
& \cdots  \xrightarrow{S}  HC_{n-1}(R) \xrightarrow{B} HH_{n}(R) \xrightarrow{I} HC_{n}(R) \xrightarrow{S} HC_{n-2}(R)  \\
& \xrightarrow{B} HH_{n-1}(R) \xrightarrow{I} HC_{n-1}(R) \to \cdots. \notag
\end{align} 

\end{theorem}

The sequence \eqref{eq:SBI} is called Connes' Periodicity Exact Sequence or SBI sequence. It naturally connects the cyclic homology $HC_{n}(R)$ with the Hochschild homology $HH_{n}(R)$.

For an associative commutative unital $k$-algebra $R$, both Hochschild homology $HH_{n}(R)$ and cyclic homology $HC_{n}(R)$ are equipped with lambda operations $\lambda^{m}$ and Adams operations $\psi^{m}$, see section 4.5 of \cite{Loday} and section 9.4.3 of \cite{W2} for details. For each integer $n \geq 1$, the group $HH_{n}(R)$ decomposes into a direct sum of eigenspaces
\begin{equation*}
HH_{n}(R)=HH^{(1)}_{n}(R) \oplus \cdots \oplus HH^{(l)}_{n}(R) \oplus \cdots \oplus HH^{(n)}_{n}(R),
\end{equation*}where $HH^{(l)}_{n}(R)$ is the eigenspace of $\psi^{m}=m^{l+1}$. 

There is a similar decomposition of cyclic homology
\begin{equation*} \label{eq:cyclic-decomp}
HC_{n}(R)=HC^{(1)}_{n}(R) \oplus \cdots \oplus HC^{(l)}_{n}(R) \oplus \cdots \oplus HC^{(n)}_{n}(R),
\end{equation*}
For $n=0$, we have the following
\begin{equation} \label{eq:hchhzeroR}
HH_{0}(R)=HH^{(0)}_{0}(R)=R, \ HC_{0}(R)=HC^{(0)}_{0}(R)=R.
\end{equation}

\begin{example} [cf. Th 4.6.10 of \cite{Loday}] \label{ex:HC-k-eigen}
The eigenspaces of cyclic homology of $k$ is known as 
\begin{align*}
\begin{cases}
 \begin{CD}
   & HC^{(l)}_{2n}(k)=k, \ \mathrm{when} \ l = n, \\
  &  HC^{(l)}_{2n}(k)=0, \ \mathrm{when} \ l \neq n.
 \end{CD}
\end{cases}
\end{align*}

\end{example}

\begin{lemma} [cf. Ex 9.4.4 and Corollary 9.8.16 of \cite{W2}]\label{lemma: l-l-omega}
With the notation as above, there are isomorphisms
\[
HH^{(n)}_{n}(R) \cong \Omega^{n}_{R /k}, \ HC^{(n)}_{n}(R) \cong \dfrac{\Omega^{n}_{R/k}}{d\Omega^{n-1}_{R/k}}.
\]
\end{lemma}

Adams operations $\psi^{m}$ decompose the SBI sequence \eqref{eq:SBI} into the sum (over $l$) of the following exact sequences
\begin{align} \label{eq:SBIeigen}
& \cdots  \xrightarrow{S}  HC^{(l-1)}_{n-1}(R)  \xrightarrow{B} HH^{(l)}_{n}(R) \xrightarrow{I} HC^{(l)}_{n}(R) \xrightarrow{S} HC^{(l-1)}_{n-2}(R)  \\
& \xrightarrow{B} HH^{(l)}_{n-1}(R) \xrightarrow{I} HC^{(l)}_{n-1}(R) \to \cdots. \notag
\end{align}

For $A \in Art_{k}$, the augmentation $A \to k$ induces a morphism 
\[
HC_{n-1}(R\otimes_{k}A)  \xrightarrow{A \to k} HC_{n-1}(R\otimes_{k}k) \xrightarrow{\cong} HC_{n-1}(R),
\]whose kernel is the relative cyclic homology,  denoted $\overline{HC}_{n-1}(R\otimes_{k}A)$. Similarly, there is a morphism of eigenspaces
\[
HC^{(l-1)}_{n-1}(R\otimes_{k}A) \xrightarrow{A \to k} HC^{(l-1)}_{n-1}(R\otimes_{k}k) \xrightarrow{\cong} HC^{(l-1)}_{n-1}(R),
\]whose kernel is denoted by $\overline{HC}^{(l-1)}_{n-1}(R\otimes_{k}A)$. 

 By \cite{Soule}, Adams operations $\psi^{m}$ exist on algebraic K-groups $K_{n}(R\otimes_{k}A)$ and $K_{n}(R)$.  We can define the relative K-group $\overline{K}_{n}(R\otimes_{k}A)$ and its eigenspaces $\overline{K}_{n}^{(l)}(R\otimes_{k}A)$ (of $\psi^{m}=m^{l}$) in a similar way as defining $\overline{HC}_{n-1}(R\otimes_{k}A)$ and $\overline{HC}^{(l-1)}_{n-1}(R\otimes_{k}A)$.

\begin{theorem} [\cite{Cath,Good}] \label{theorem: GoodwillieCathelineau}
With the notation as above, the relative Chern character induces the following isomorphisms $\mathrm{(}$where the cyclic homology is defined over $\mathbb{Q} \mathrm{)}$,
\begin{align*}
\begin{cases}
 \overline{K}_{n}(R\otimes_{k}A)_{\mathbb{Q}} \xrightarrow{\cong} \overline{HC}_{n-1}(R\otimes_{k}A), \\
\overline{K}_{n}^{(l)}(R\otimes_{k}A)_{\mathbb{Q}}  \xrightarrow{\cong} \overline{HC}_{n-1}^{(l-1)}(R\otimes_{k}A).
\end{cases}
\end{align*}

\end{theorem}

This theorem is very useful for computing relative K-groups, and we have used it to study deformations of Chow groups in \cite{DHY,Y2,Y5}. Corti$\mathrm{\tilde{n}}$as-Haesemeyer-Weibel \cite{CHW} generalized this theorem to the space level. 

In the remainder of this section, $X$ is a smooth projective variety over an arbitrary field $k$ of characteristic zero. For any $A \in Art_{k}$, we write $X_{A}=X \otimes_{\mathrm{Spec}(k)} \mathrm{Spec}(A)$.

\begin{definition}  \label{d:HH-HC-sheaf}
Let  $HH_{n}(O_{X_{A}})$ and $HC_{n}(O_{X_{A}})$ be the Hochschild  and cyclic homology sheaves associated to the following presheaves respectively 
\[
U \to HH_{n}(O_{X_{A}}(U)), \ U \to HC_{n}(O_{X_{A}}(U)), 
\]where $U=\mathrm{Spec}(R) \subset X$ is open affine, $HH_{n}(O_{X_{A}}(U))=HH_{n}(R\otimes_{k} A)$ and $HC_{n}(O_{X_{A}}(U))=HC_{n}(R\otimes_{k} A)$ are defined in \eqref{eq:hhalg} and \eqref{eq:hcalg} respectively.

\end{definition}

Adams operations $\psi^{m}$ act on the sheaves  $HC_{n}(O_{X_{A}})$ and $HH_{n}(O_{X_{A}})$, whose eigenspaces (of $\psi^{m}=m^{l+1}$) are denoted by $HC^{(l)}_{n}(O_{X_{A}})$ and $HH^{(l)}_{n}(O_{X_{A}})$ respectively. When $A=k$, we write $HC^{(l)}_{n}(O_{X})= HC^{(l)}_{n}(O_{X_{k}})$ and $HH^{(l)}_{n}(O_{X})= HH^{(l)}_{n}(O_{X_{k}})$.

\begin{definition}  \label{d:rel-HH-HC-eigen}
 The relative cyclic and Hochschild homology sheaves $\overline{HC}^{(l)}_{n}(O_{X_{A}})$ and $\overline{HH}^{(l)}_{n}(O_{X_{A}})$ are defined to be kernels of the following maps respectively (induced by $A \to k$) 
\[
 HC^{(l)}_{n}(O_{X_{A}})  \xrightarrow{A \to k} HC^{(l)}_{n}(O_{X}), \ HH^{(l)}_{n}(O_{X_{A}}) \xrightarrow{A \to k} HH^{(l)}_{n}(O_{X}).
\]
\end{definition}

It is clear that the following short exact sequence splits,
\begin{equation} \label{eq:ses-hh}
 0 \to \overline{HH}^{(l)}_{n}(O_{X_{A}}) \to HH^{(l)}_{n}(O_{X_{A}}) \xrightarrow{A \to k} HH^{(l)}_{n}(O_{X}) \to 0.
\end{equation}

Now we use the relative cyclic homology sheaf to compute the formal completion $\widehat{CH}^{p}(A)$, which is defined in \eqref{eq:formalcom}. 
\begin{lemma} [cf. Lemma 3.5 of \cite{Y6}] \label{l:iso2}
With the notation as above, there is an isomorphism
\begin{equation*}
\widehat{CH}^{p}(A) \cong H^{p}(X, \overline{HC}^{(p-1)}_{p-1}(O_{X_{A}})).
\end{equation*}

\end{lemma}

We briefly outline the proof for the convenience of the reader. By \cite{Soule}, the Quillen K-theory sheaf $K_{p}(O_{X_{A}})$ carries Adams operations $\psi^{m}$, whose eigenspaces (of $\psi^{m}=m^{l}$) are $K^{(l)}_{p}(O_{X_{A}})$, where $1 \leq l \leq p$.

Let $\overline{K}_{p}^{M}(O_{X_{A}})$ and $\overline{K}_{p}^{(p)}(O_{X_{A}})$  be the kernels of the following maps respectively 
\[
K_{p}^{M}(O_{X_{A}}) \xrightarrow{A \to k} K_{p}^{M}(O_{X}), \ K_{p}^{(p)}(O_{X_{A}}) \xrightarrow{A \to k} K_{p}^{(p)}(O_{X}),
\]where $K_{p}^{M}(O_{X_{A}})$ is the Milnor K-theory sheaf.

By a sheaf version of Theorem 5 on page 526 of \cite{Soule}, there exists isomorphisms $K_{p}^{M}(O_{X_{A}})_{\mathbb{Q}} \xrightarrow{\cong} K_{p}^{(p)}(O_{X_{A}})_{\mathbb{Q}}$ and $K_{p}^{M}(O_{X})_{\mathbb{Q}} \xrightarrow{\cong} K_{p}^{(p)}(O_{X})_{\mathbb{Q}}$. This yields that 
\begin{equation} \label{eq:Kseq}
\overline{K}_{p}^{M}(O_{X_{A}})_{\mathbb{Q}}  \xrightarrow{\cong} \overline{K}_{p}^{(p)}(O_{X_{A}})_{\mathbb{Q}} \cong \overline{HC}^{(p-1)}_{p-1}(O_{X_{A}}),
\end{equation}where the second isomorphism is a sheaf version of Theorem \ref{theorem: GoodwillieCathelineau} (let $l=n=p$).

Since the short exact sequence 
\[
0 \to \overline{K}_{p}^{M}(O_{X_{A}}) \to K_{p}^{M}(O_{X_{A}}) \to K_{p}^{M}(O_{X}) \to 0
\]splits, it follows that
\begin{equation*} \label{eq:iso1}
\widehat{CH}^{p}(A)  \cong H^{p}(X,\overline{K}_{p}^{M}(O_{X_{A}}))_{\mathbb{Q}}.
\end{equation*}We further have
{\small
\begin{equation*} \label{eq:iso1}
\widehat{CH}^{p}(A)  \cong H^{p}(X,\overline{K}_{p}^{M}(O_{X_{A}}))_{\mathbb{Q}} \cong  H^{p}(X,\overline{K}_{p}^{(p)}(O_{X_{A}}))_{\mathbb{Q}} \cong  H^{p}(X,\overline{HC}^{(p-1)}_{p-1}(O_{X_{A}})),
\end{equation*}
}where the last two isomorphisms are from \eqref{eq:Kseq}.

\section{Deformations of Chow groups}

In this section, $X$ is a smooth projective variety over an arbitrary field $k$ of characteristic zero. For $A \in  Art_{k}$, we write $X_A=X\times_{\mathrm{Spec}(k)}\mathrm{Spec}(A)$ as before. We study deformations of the Chow group $CH^{p}(X)$ by computing its formal completion $\widehat{CH}^{p}(A)$, where $p$ is an integer such that $1 \leq p \leq \mathrm{dim}(X)$.

\begin{definition} \label{d:rep}
A functor $F$ is called pro-representable if it is isomorphic to the functor $h_{R}$, which is defined as 
\[
h_{R}(A)=Hom(R,A).
\]Here $A \in  Art_{k}$, $R$ is a complete local Noetherian $k$-algebra with the maximal ideal $m_{R}$ such that $R/m_{R}^{n} \in  Art_{k}$, where $n$ is any positive integer, and $Hom(R,A)$ denotes the set of local $k$-algebra homomorphisms from $R$ to $A$. 

\end{definition}

A useful example of a pro-representable functor $\mathbb{T}$ is given as
\begin{equation} \label{eq: pro-rep vs}
\mathbb{T}(A) = V \otimes_{k}m_{A},
\end{equation}where $V$ is a finite dimensional $k$-vector space, and $m_{A}$ is the maximal ideal of $A \in Art_{k}$. We refer to \cite{Sc} for details on pro-representable functors.

To compute the formal completion $\widehat{CH}^{p}(A)$, we begin with the case of $p=1$, which is classical. By Lemma \ref{l:iso2} (let $p=1$), there is an isomorphism
\begin{equation*}
\widehat{CH}^{1}(A) \cong H^{1}(X, \overline{HC}^{(0)}_{0}(O_{X_{A}})).
\end{equation*} 

The relative cyclic homology sheaf $\overline{HC}^{(0)}_{0}(O_{X_{A}})$ is the kernel of 
\[
HC^{(0)}_{0}(O_{X_{A}}) \xrightarrow{A \to k} HC^{(0)}_{0}(O_{X}), 
\]cf. Definition \ref{d:rel-HH-HC-eigen}. A sheaf version of \eqref{eq:hchhzeroR} says that
\[
HC^{(0)}_{0}(O_{X_{A}}) \cong O_{X_{A}}\cong O_{X} \otimes_{k}A,  \ HC^{(0)}_{0}(O_{X}) \cong O_{X}.
\]Hence, $\overline{HC}^{(0)}_{0}(O_{X_{A}})$ is isomorphic to the kernel of 
\[
 O_{X} \otimes_{k}A \xrightarrow{A \to k}O_{X},
\]which can be computed directly as 
\begin{equation} \label{eq:hc00}
\overline{HC}^{(0)}_{0}(O_{X_{A}}) \cong O_{X} \otimes_{k}m_{A}.
\end{equation}

As a consequence, we have
\begin{equation*}
\widehat{CH}^{1}(A) \cong H^{1}(X, O_{X} \otimes_{k}m_{A}) \cong H^{1}(X, O_{X}) \otimes_{k}m_{A}.
\end{equation*}

Since the functor $\widehat{CH}^{1}$ carries the form of \eqref{eq: pro-rep vs}, where $V=H^{1}(X, O_{X})$, the pro-representability of $\widehat{CH}^{1}$ is clear. To summarize, we have the following result, which should be known to experts.
\begin{lemma} \label{l:ch1}
With the notation as above, there is an isomorphism
\[
\widehat{CH}^{1}(A) \cong H^{1}(X, O_{X}) \otimes_{k}m_{A}.
\]Hence, the functor $\widehat{CH}^{1}$ is pro-representable.

\end{lemma}

Using the isomorphisms $K^{M}_{1}(O_{X_{A}}) \cong O_{X_{A}}^{*}$ and $K^{M}_{1}(O_{X}) \cong O_{X}^{*}$, we can show that $\widehat{CH}^{1}(A)$, defined in \eqref{eq:formalcom}, is isomorphic to the kernel of the morphism
\[
H^{1}(X, O_{X_{A}}^{*}) \xrightarrow{A \to k} H^{1}(X, O_{X}^{*}),
\]which is $H^{1}(X, O_{X}) \otimes_{k}m_{A}$. This provides an alternative proof of Lemma \ref{l:ch1}.

In the remainder of this section, we compute the functor $\widehat{CH}^{p}$, where $p \geq 2$. Since $ p \leq \mathrm{dim}(X)$, this forces $\mathrm{dim}(X) \geq 2$. We adopt the following setting.

\begin{setting} [cf. Setting 4.1 of \cite{Y6}]  \label{s:set}
Let $X$ be a smooth projective variety of dimension \(d \geq 2\) over a field $k$ of characteristic zero. Let $p$ be an integer such that $ 2 \leq p \leq d$. We assume that $X$ satisfies the following condition:
\begin{equation} \label{eq: condi-vanish}
H^{p}(X, \Omega^{i}_{X/k})= H^{p+1}(X, \Omega^{i}_{X/k})= \cdots =H^{2p-1-i}(X, \Omega^{i}_{X/k})=0,
\end{equation}where $i$ is an integer such that $0 \leq i \leq p-2$.

\end{setting}

The condition \eqref{eq: condi-vanish}, which has been recalled as \eqref{eq: condi-vanish1} in Theorem \ref{t:b-m-yang}, was firstly introduced in Setting 4.1 of \cite{Y6}. It can be alternatively written as
 {\footnotesize
\begin{align}\label{eq: condi-vanish-explicit}
\begin{cases}
   i=0: \ H^{p}(X, O_{X}) = H^{p+1}(X, O_{X})= \cdots  = H^{2p-1}(X,O_{X})=0, \\
   i=1: \ H^{p}(X,\Omega^{1}_{X/k}) = H^{p+1}(X,\Omega^{1}_{X/k})= \cdots =H^{2p-2}(X,\Omega^{1}_{X/k})=0, \\
    \ \ \ \ \ \   \vdots \\
i=p-2: \ H^{p}(X,\Omega^{p-2}_{X/k}) = H^{p+1}(X,\Omega^{p-2}_{X/k})=0.
\end{cases}
\end{align}
}We refer to section 4 of \cite{Y6} for more discussions on this condition.

For any $A \in Art_{k}$, we write $X_{A}=X \otimes_{\mathrm{Spec}(k)} \mathrm{Spec}(A)$ as before. Since $X$ is smooth over $k$, by base change, $X_{A}$ is smooth over $A$. There is a short exact sequence of sheaves
\[
0 \to \Omega^{1}_{A/ \mathbb{Q}}\otimes_{A}O_{X_{A}} \to \Omega^{1}_{X_{A}/ \mathbb{Q}} \to \Omega^{1}_{X_{A}/ A} \to 0.
\]This yields a filtration on $\Omega^{i}_{X_{A}/ \mathbb{Q}}$, where $i$ is the integer in Setting \ref{s:set}. To be precise, let $F^{j}$ be the image of the map 
\[
\Omega^{j}_{A/ \mathbb{Q}} \otimes_{A} \Omega^{i-j}_{X_{A}/ \mathbb{Q}} \to \Omega^{i}_{X_{A}/ \mathbb{Q}},
\]where $j$ is an integer satisfying that $0 \leq j \leq i$, there is a filtration on $\Omega^{i}_{X_{A}/ \mathbb{Q}}$ given by
\begin{equation*}
\Omega^{i}_{X_{A}/ \mathbb{Q}}=F^{0} \supset F^{1} \supset \cdots \supset F^{i} \supset F^{i+1}=0.
\end{equation*}

The associated graded piece is 
\[
Gr^{j}\Omega^{i}_{X_{A}/ \mathbb{Q}}=F^{j}/F^{j+1}=\Omega^{j}_{A/ \mathbb{Q}} \otimes_{A} \Omega^{i-j}_{X_{A}/ A}.
\]We use the isomorphism $\Omega^{i-j}_{X_{A}/ A} \cong \Omega^{i-j}_{X/ k}\otimes_{k} A$ to deduce that
\begin{align}  \label{eq:gradeiso1}
Gr^{j}\Omega^{i}_{X_{A}/ \mathbb{Q}} \cong \Omega^{j}_{A/ \mathbb{Q}} \otimes_{A} (\Omega^{i-j}_{X/ k}\otimes_{k} A) \cong \Omega^{i-j}_{X/ k} \otimes_{k}\Omega^{j}_{A/ \mathbb{Q}}.
\end{align}

\begin{lemma}  \label{l:extvanish}
In Setting \ref{s:set}, with the notation as above, we have the following
\begin{equation*}
H^{p}(X, F^{j})= \cdots =H^{2p-1-(i-j)}(X, F^{j})=0.
\end{equation*}In particular, when $j=0$ $\mathrm{(}F^{j}=F^{0}=\Omega^{i}_{X_{A}/ \mathbb{Q}}$ in this case$\mathrm{)}$, we have
\begin{equation*}
H^{p}(X, \Omega^{i}_{X_{A}/ \mathbb{Q}})= \cdots =H^{2p-1-i}(X, \Omega^{i}_{X_{A}/ \mathbb{Q}})=0.
\end{equation*}
\end{lemma}

\begin{proof}
This can be proved by induction. We start with the case of $j=i$. When $j=i$, $2p-1-(i-j)=2p-1$. Note that $F^{i}=Gr^{i}\Omega^{i}_{X_{A}/ \mathbb{Q}}=O_{X} \otimes_{k} \Omega^{i}_{A/ \mathbb{Q}}$, cf. \eqref{eq:gradeiso1}, the condition \eqref{eq: condi-vanish} yields that
{\small
\begin{equation*}  \label{eq:zerocohfil}
H^{q}(X, F^{i})\cong H^{q}(X, O_{X})\otimes_{k}\Omega^{i}_{A/ \mathbb{Q}}=0,  \ \mathrm{where} \ q=p, \ p+1, \ \cdots, \ 2p-1.
\end{equation*}
}Hence, the case of $j=i$ is verified.

Now, we suppose that the case of $j+1$ has been proved, where $0 \leq j \leq i-1$, i.e., we have
\begin{equation}  \label{eq:zerocohj}
H^{q}(X, F^{j+1})=0,  \ \mathrm{where} \ q=p, \ p+1, \ \cdots, \ 2p-1-(i-(j+1)).
\end{equation}

We next compute the cohomology $H^{q}(X, F^{j})$. For that purpose, we consider the following short exact sequence of sheaves 
\begin{equation*}  \label{eq:sesfil}
0 \to F^{j+1} \to F^{j} \to Gr^{j}\Omega^{i}_{X_{A}/ \mathbb{Q}} \to 0,
\end{equation*}whose associated long exact sequence of sheaf cohomology groups is
{\footnotesize
\begin{equation} \label{eq:longseq}
\cdots \to H^{q}(X, F^{j+1}) \to H^{q}(X, F^{j}) \to H^{q}(X, Gr^{j}\Omega^{i}_{X_{A}/ \mathbb{Q}}) \to \cdots.
\end{equation}
}

Since $0 \leq i-j \leq i \leq p-2$, for $q=p$, $p+1$, $\cdots$, $2p-1-(i-j)$, the condition \eqref{eq: condi-vanish} says that
\begin{equation*} \label{eq:zerocoh}
H^{q}(X,\Omega^{i-j}_{X/ k})=0.
\end{equation*}Combining it with the isomorphism \eqref{eq:gradeiso1}, we have
 \begin{align}\label{eq:zerocoh2}
  H^{q}(X, Gr^{j}\Omega^{i}_{X_{A}/ \mathbb{Q}}) & \cong  H^{q}(X,\Omega^{i-j}_{X/ k} \otimes_{k}\Omega^{j}_{A/ \mathbb{Q}})  \\
& \cong  H^{q}(X,\Omega^{i-j}_{X/ k}) \otimes_{k}\Omega^{j}_{A/ \mathbb{Q}}=0. \notag
\end{align}

Substituting \eqref{eq:zerocohj} and \eqref{eq:zerocoh2} into this sequence \eqref{eq:longseq}, we see that
\begin{equation}  \label{eq:zerocohfil6}
H^{q}(X, F^{j})=0,  \ \mathrm{where} \ q=p, \ p+1, \ \cdots, \ 2p-1-(i-j).
\end{equation}This finishes the induction. 

In particular, when $j=0$, $F^{j}=F^{0}=\Omega^{i}_{X_{A}/ \mathbb{Q}}$. Note that $i-j=i$ in this case, and \eqref{eq:zerocohfil6} becomes as
\begin{equation*}
H^{p}(X, \Omega^{i}_{X_{A}/ \mathbb{Q}})= \cdots =H^{2p-1-i}(X, \Omega^{i}_{X_{A}/ \mathbb{Q}})=0.
\end{equation*}

\end{proof}

By Lemma \ref{lemma: l-l-omega} (let $n=i$), there is an isomorphism
\[
 HH^{(i)}_{i}(O_{X_{A}}) \cong \Omega^{i}_{X_A/ \mathbb{Q}},
\]where the Hochschild homology is defined over $\mathbb{Q}$. Hence, Lemma \ref{l:extvanish} immediately yields that
\begin{equation} \label{eq:HH1}
H^{p}(X, HH^{(i)}_{i}(O_{X_{A}}))= \cdots =H^{2p-1-i}(X, HH^{(i)}_{i}(O_{X_{A}}))=0.
\end{equation}

When $A=k$, the above \eqref{eq:HH1} specializes into 
\begin{equation} \label{eq:HH2}
H^{p}(X, HH^{(i)}_{i}(O_{X}))= \cdots =H^{2p-1-i}(X, HH^{(i)}_{i}(O_{X}))=0.
\end{equation}

Let $n=l=i$ in the split short exact sequence \eqref{eq:ses-hh}. It yields a short exact sequence of cohomology groups
{\footnotesize
\[
0 \to H^{j}(X,\overline{HH}^{(i)}_{i}(O_{X_{A}})) \to H^{j}(X,HH^{(i)}_{i}(O_{X_{A}})) \to  H^{j}(X,HH^{(i)}_{i}(O_{X})) \to 0,
\]
}where $j$ is a non-negative integer. As a consequence, we use  \eqref{eq:HH1} to deduce the following result.

\begin{lemma} \label{l:cohrelhhzero}
In Setting \ref{s:set}, where $i=0$, $1$, $\cdots$, $p-2$, we have
\begin{equation} \label{eq:rel-HH-coh0}
H^{q}(X,\overline{HH}^{(i)}_{i}(O_{X_{A}}))=0,\ \mathrm{where} \ q=p, \ p+1, \ \cdots, \ 2p-1-i.
\end{equation}

\end{lemma}
We can alternatively write the isomorphism \eqref{eq:rel-HH-coh0} as 
 {\footnotesize
\begin{align*}\label{eq: HH-vanish-explicit}
\begin{cases}
   i=0: \ H^{p}(X, \overline{HH}^{(0)}_{0}(O_{X_{A}})) = \cdots  = H^{2p-1}(X,\overline{HH}^{(0)}_{0}(O_{X_{A}}))=0, \\
   i=1: \ H^{p}(X,\overline{HH}^{(1)}_{1}(O_{X_{A}})) =  \cdots =H^{2p-2}(X,\overline{HH}^{(1)}_{1}(O_{X_{A}}))=0, \\
    \ \ \ \ \ \   \vdots \\
i=p-2: \ H^{p}(X,\overline{HH}^{(p-2)}_{p-2}(O_{X_{A}})) = H^{p+1}(X,\overline{HH}^{(p-2)}_{p-2}(O_{X_{A}}))=0.
\end{cases}
\end{align*}
}

We next compute the case of $p-1$, i.e., we compute the cohomology groups $H^{p}(X, HH^{(p-1)}_{p-1}(O_{X_{A}}))$ and $H^{p}(X, \overline{HH}^{(p-1)}_{p-1}(O_{X_{A}}))$.

\begin{lemma} 
In Setting \ref{s:set}, there is an isomorphism
\begin{equation} \label{eq:cohomegaA}
H^{p}(X,\Omega^{p-1}_{X_{A}/ \mathbb{Q}}) \cong H^{p}(X,\Omega^{p-1}_{X/ k })\otimes_{k} A.
\end{equation}

\end{lemma}

This can be proved similarly to Lemma  \ref{l:extvanish}. We write it out for readers' convenience.

\begin{proof}
Let $F^{j}$ be the image of the map 
\[
\Omega^{j}_{A/ \mathbb{Q}} \otimes_{A} \Omega^{p-1-j}_{X_{A}/ \mathbb{Q}} \to \Omega^{p-1}_{X_{A}/ \mathbb{Q}},
\]where $j$ is an integer satisfying that $0 \leq j \leq p-1$. There is a filtration on $\Omega^{p-1}_{X_{A} / \mathbb{Q}}$ given by
\begin{equation*}
\Omega^{p-1}_{X_{A}/ \mathbb{Q}}=F^{0} \supset F^{1} \supset \cdots \supset F^{p-1} \supset F^{p}=0,
\end{equation*}
whose associated graded piece is 
\begin{align} \label{eq:gradeiso}
Gr^{j}\Omega^{p-1}_{X_{A} / \mathbb{Q}}& =F^{j}/F^{j+1}=\Omega^{j}_{A/ \mathbb{Q}} \otimes_{A} \Omega^{p-1-j}_{X_{A}/ A}\\ \notag
& \cong \Omega^{j}_{A/ \mathbb{Q}} \otimes_{A} (\Omega^{p-1-j}_{X/ k} \otimes_{k} A) \\ \notag
& \cong  \Omega^{p-1-j}_{X/ k} \otimes_{k} \Omega^{j}_{A/ \mathbb{Q}}.
\end{align}

We first consider the case of $j$ satisfying that $1 \leq j \leq p-1$, and leave the case of $j=0$ to the end of the proof.

When $j$ satisfies that $1 \leq j \leq p-1$, $0 \leq p-1-j \leq p-2$. By letting $i=p-1-j$ in the condition \eqref{eq: condi-vanish}, we have
\begin{equation*}
H^{q}(X, \Omega^{p-1-j}_{X/ k})=0,
\end{equation*}where $q=p$, $p+1$, $\cdots$, $p+j$ ($=2p-1-i=2p-1-(p-1-j)$).

We use the isomorphism \eqref{eq:gradeiso} to deduce that
\begin{equation} \label{eq:gradecohzero}
H^{q}(X, Gr^{j}\Omega^{p-1}_{X_{A} / \mathbb{Q}}) \cong H^{q}(X, \Omega^{p-1-j}_{X/ k}) \otimes_{k} \Omega^{j}_{A/ \mathbb{Q}}=0,
\end{equation}where $q=p,\ \cdots, \ p+j$.

When $j=p-1$, $F^{p-1}=Gr^{p-1}\Omega^{p-1}_{X_{A}/ \mathbb{Q}} \cong O_{X}\otimes_{k} \Omega^{p-1}_{A/ \mathbb{Q}}$. The condition \eqref{eq: condi-vanish} yields that 
{\small
\begin{equation} \label{eq:filcohzero}
H^{q}(X, F^{p-1})\cong H^{q}(X,O_{X}\otimes_{k} \Omega^{p-1}_{A/ \mathbb{Q}})\cong H^{q}(X,O_{X})\otimes_{k} \Omega^{p-1}_{A/ \mathbb{Q}}=0, 
\end{equation}
}where $q=p$, $\cdots$, $2p-1 (=p+j=p+(p-1))$.

When $j=p-2$, there is a short exact sequence of sheaves
\begin{equation*} \label{eq:sesfil2}
0 \to F^{p-1} \to F^{p-2} \to Gr^{p-2}\Omega^{p-1}_{X_{A}/ \mathbb{Q}} \to 0,
\end{equation*}whose associated long exact sequence of sheaf cohomology is
{\footnotesize
\begin{equation} \label{eq:lesfil2}
\cdots \to H^{q}(X, F^{p-1}) \to H^{q}(X, F^{p-2}) \to H^{q}(X, Gr^{p-2}\Omega^{p-1}_{X_{A}/ \mathbb{Q}}) \to \cdots.
\end{equation}
}

Let $j=p-2$ in \eqref{eq:gradecohzero}, we have
\begin{equation} \label{eq:zerocohgrade22}
H^{q}(X, Gr^{p-2}\Omega^{p-1}_{X_{A}/ \mathbb{Q}})=0, 
\end{equation}where $q=p$, $\cdots$, $2p-2 (=p+j=p+(p-2))$.

Substituting \eqref{eq:filcohzero} and \eqref{eq:zerocohgrade22} into the sequence \eqref{eq:lesfil2}, we see that 
\begin{equation*}
H^{q}(X,F^{p-2})=0, \ \mathrm{where}\ q=p,\ \cdots, \ 2p-2.
\end{equation*}

There are short exact sequences of sheaves
\begin{equation*}
0 \to F^{j+1} \to F^{j} \to Gr^{j}\Omega^{p-1}_{X_{A}/ \mathbb{Q}} \to 0,
\end{equation*}where $j=1, \cdots, p-3$. By continuing the above procedure, we can prove the following
\begin{equation*} \label{eq:zerocohgrade2}
H^{q}(X,F^{j})=0, \  \mathrm{where}\ q=p,\cdots, \ p+j.
\end{equation*}In particular, when $j=1$, we have
\begin{equation} \label{eq:cohfil1zero}
H^{p}(X, F^{1})=H^{p+1}(X,F^{1})=0.
\end{equation}

Now, we handle the case of $F^{0}$. There is a short exact sequence of sheaves
\begin{equation*} \label{eq:sesfil0}
0 \to F^{1} \to F^{0} \to Gr^{0}\Omega^{p-1}_{X_{A}/ \mathbb{Q}} \to 0,
\end{equation*}whose associated long exact sequence of sheaf cohomologies is
{\footnotesize
\begin{equation*} \label{eq:les01}
\cdots \to H^{p}(X, F^{1}) \to H^{p}(X, F^{0}) \to H^{p}(X,Gr^{0}\Omega^{p-1}_{X/ \mathbb{Q}}) \to H^{p+1}(X,F^{1}) \to \cdots.
\end{equation*}
}Substituting \eqref{eq:cohfil1zero} into this sequence, we have 
\begin{equation*}
H^{p}(X, F^{0}) \cong H^{p}(X, Gr^{0}\Omega^{p-1}_{X/ \mathbb{Q}}) \cong H^{p}(X, \Omega^{p-1}_{X/ k})\otimes_{k} A,
\end{equation*}where the second isomorphism is induced from  \eqref{eq:gradeiso} (let $j=0$).

Since $F^{0}=\Omega^{p-1}_{X_{A}/ \mathbb{Q}}$, this completes the proof.

\end{proof}

When $A=k$, the isomorphism \eqref{eq:cohomegaA} carries the form
\begin{equation} \label{eq:cohomega}
H^{p}(X, \Omega^{p-1}_{X/ \mathbb{Q}}) \cong H^{p}(X, \Omega^{p-1}_{X/ k})\otimes_{k} k \cong H^{p}(X, \Omega^{p-1}_{X/ k}).   
\end{equation}

By Lemma \ref{lemma: l-l-omega} (let $n=p-1$), we rewrite \eqref{eq:cohomegaA} and \eqref{eq:cohomega} as 
\begin{align}\label{eq:cohomegahh}
\begin{cases}
H^{p}(X, HH^{(p-1)}_{p-1}(O_{X_{A}}))=H^{p}(X, \Omega^{p-1}_{X/ k})\otimes_{k} A  \\
H^{p}(X, HH^{(p-1)}_{p-1}(O_{X}))=H^{p}(X, \Omega^{p-1}_{X/ k}),
\end{cases}
\end{align}where the Hochschild homologies are defined over $\mathbb{Q}$.

Let $n=l=p-1$ in the short exact sequence \eqref{eq:ses-hh}. Since it is split, it yields a short exact sequence of cohomology groups
{\footnotesize
\[
0 \to H^{p}(X,\overline{HH}^{(p-1)}_{p-1}(O_{X_{A}})) \to H^{p}(X,HH^{(p-1)}_{p-1}(O_{X_{A}})) \xrightarrow{A \to k}  H^{p}(X,HH^{(p-1)}_{p-1}(O_{X})) \to 0.
\]
}Substituting \eqref{eq:cohomegahh} into it, we rewrite it as 
\[
0 \to H^{p}(X,\overline{HH}^{(p-1)}_{p-1}(O_{X_{A}})) \to H^{p}(X, \Omega^{p-1}_{X/ k})\otimes_{k} A \xrightarrow{A \to k}  H^{p}(X, \Omega^{p-1}_{X/ k}) \to 0.
\]Note that $A$ decomposes into a direct sum of $k$-vector spaces
\[
A=k \oplus m_{A},
\]where $m_{A}$ is the maximal ideal of $A$. It follows that
\begin{align*}
H^{p}(X, \Omega^{p-1}_{X/ k})\otimes_{k} A & \cong H^{p}(X, \Omega^{p-1}_{X/ k})\otimes_{k} (k \oplus m_{A}) \\
& \cong (H^{p}(X, \Omega^{p-1}_{X/ k})\otimes_{k} k) \oplus (H^{p}(X, \Omega^{p-1}_{X/ k})\otimes_{k}m_{A}) \\
& \cong H^{p}(X, \Omega^{p-1}_{X/ k})\oplus (H^{p}(X, \Omega^{p-1}_{X/ k})\otimes_{k}m_{A}).
\end{align*}
Consequently, we compute $H^{p}(X,\overline{HH}^{(p-1)}_{p-1}(O_{X_{A}}))$ as following.
\begin{lemma} \label{l:cohrelhh}
In Setting \ref{s:set}, there is an isomorphism
\[
H^{p}(X, \overline{HH}^{(p-1)}_{p-1}(O_{X_{A}})) \cong H^{p}(X,\Omega^{p-1}_{X/ k})\otimes_{k}{m_{A}}.
\]
\end{lemma}

From now on, for $A \in Art_{k}$, we further assume that $A$ is also a graded $k$-algebra $A=\oplus_{l \geq 0}A_{l}$ such that $A_{0}=k$.

Let $U=\mathrm{Spec}(R) \subset X$ be open affine, $R \otimes_{k}A= \oplus_{l \geq 0}(R \otimes_{k}A_{l})$ is a graded $k$-algebra with $R \otimes_{k}A_{0}=R \otimes_{k}k=R$. Since $\mathbb{Q} \subset k$, $R \otimes_{k}A$ can be also considered as a graded $\mathbb{Q}$-algebra.

Let $l=n$ in the SBI sequence \eqref{eq:SBIeigen}, where Hochschild and cyclic homology are defined over $\mathbb{Q}$. By Goodwillie \cite{GoodW}, it breaks into a short exact sequence 
{\footnotesize
\begin{equation} \label{eq:SBI-alg}
0 \to \overline{HC}^{(n-1)}_{n-1}(R \otimes_{k}A) \xrightarrow{\mathrm{B}} \overline{HH}^{(n)}_{n}(R \otimes_{k}A)  \xrightarrow{\mathrm{I}} \overline{HC}^{(n)}_{n}(R \otimes_{k}A) \to 0,
\end{equation}
}where $n$ is any positive integer, $\overline{HC}^{(n-1)}_{n-1}(R \otimes_{k}A)$ is defined to be the kernel of 
\[
HC^{(n-1)}_{n-1}(R \otimes_{k}A) \xrightarrow{A \to k} HC^{(n-1)}_{n-1}(R),
\]and $\overline{HH}^{(n)}_{n}(R \otimes_{k}A)$ and $\overline{HC}^{(n)}_{n}(R \otimes_{k}A)$ are defined similarly.

This sequence \eqref{eq:SBI-alg} is useful for computing cyclic homology and K-theory, for example, see Geller, Reid and Weibel \cite{GRW, GW}, and the author \cite{Y5}. 

The following short exact sequence
{\footnotesize
\begin{equation*} \label{eq:SBI-sheaf}
0 \to \overline{HC}^{(n-1)}_{n-1}(O_{X_{A}}) \xrightarrow{\mathrm{B}} \overline{HH}^{(n)}_{n}(O_{X_{A}})  \xrightarrow{\mathrm{I}} \overline{HC}^{(n)}_{n}(O_{X_{A}}) \to 0,
\end{equation*}
}is a sheaf version of \eqref{eq:SBI-alg}, where $ \overline{HH}^{(n)}_{n}(O_{X_{A}})$ etc are the relative Hochschild  and cyclic homology sheaves, cf. Definition \ref{d:rel-HH-HC-eigen}. Its associated long exact sequence of sheaf cohomology groups is
\begin{align} \label{l:lexhhhc}
& \cdots \to H^{q}(X, \overline{HC}^{(n-1)}_{n-1}(O_{X_{A}})) \to H^{q}(X,\overline{HH}^{(n)}_{n}(O_{X_{A}}))  \\
& \to H^{q}(X,\overline{HC}^{(n)}_{n}(O_{X_{A}})) \to H^{q+1}(X,\overline{HC}^{(n-1)}_{n-1}(O_{X_{A}})) \to  \cdots,  \notag
\end{align}where $q$ is a non-negative integer.

When $n=1$, $n-1=0$, the sequence \eqref{l:lexhhhc} has the form
\begin{align} \label{eq:lexhhhc1}
& \cdots \to H^{q}(X,\overline{HC}^{(0)}_{0}(O_{X_{A}})) \to H^{q}(X,\overline{HH}^{(1)}_{1}(O_{X_{A}}))  \\
& \to H^{q}(X,\overline{HC}^{(1)}_{1}(O_{X_{A}})) \to H^{q+1}(X,\overline{HC}^{(0)}_{0}(O_{X_{A}})) \to  \cdots.  \notag
\end{align} 

There is an isomorphism (cf. \eqref{eq:hc00})
\[
\overline{HC}^{(0)}_{0}(O_{X_{A}})  \cong O_{X} \otimes_{k} m_{A}.
\]Hence, the condition \eqref{eq: condi-vanish} yields that, for $q= p$, $p+1$, $\cdots$, $2p-1$, 
\begin{equation} \label{eq:cohrelhc}
H^{q}(X,\overline{HC}^{(0)}_{0}(O_{X_{A}})) \cong H^{q}(X,O_{X}) \otimes_{k} m_{A}=0.
\end{equation}

Substituting \eqref{eq:cohrelhc} into the sequence \eqref{eq:lexhhhc1}, for $q=p, p+1, \cdots, 2p-2$, we see that
\begin{align*}
H^{q}(X,\overline{HC}^{(1)}_{1}(O_{X_{A}})) \cong H^{q}(X,\overline{HH}^{(1)}_{1}(O_{X_{A}}))=0,
\end{align*}where the second identity is from  Lemma \ref{l:cohrelhhzero} (let $i=1$). That is, we have
\begin{align} \label{eq:cohrelhc1}
H^{q}(X,\overline{HC}^{(1)}_{1}(O_{X_{A}}))=0, \ \mathrm{where} \ q=p, \ p+1, \ \cdots, \ 2p-2.
\end{align} 

Let $n=2$ in the sequence \eqref{l:lexhhhc}, it has the form
\begin{align*} 
& \cdots \to H^{q}(X, \overline{HC}^{(1)}_{1}(O_{X_{A}})) \to H^{q}(X,\overline{HH}^{(2)}_{2}(O_{X_{A}}))  \\
& \to H^{q}(X, \overline{HC}^{(2)}_{2}(O_{X_{A}})) \to H^{q+1}(X,\overline{HC}^{(1)}_{1}(O_{X_{A}})) \to  \cdots.  \notag
\end{align*}

Substituting \eqref{eq:cohrelhc1} into this sequence, for $q=p, p+1, \cdots, 2p-3$, we see that
\begin{align*}
H^{q}(X, \overline{HC}^{(2)}_{2}(O_{X_{A}})) \cong H^{q}(X, \overline{HH}^{(2)}_{2}(O_{X_{A}}))=0,
\end{align*}where the second identity is from  Lemma \ref{l:cohrelhhzero} (let $i=2$). That is, we have
\begin{align*}
H^{q}(X, \overline{HC}^{(2)}_{2}(O_{X_{A}}))=0, \ \mathrm{where} \ q=p, \ p+1, \ \cdots, \ 2p-3.
\end{align*}

We continue this procedure by letting $n=3,\ 4, \ \cdots, \ p-2$ in the sequence \eqref{l:lexhhhc}, and can prove that, for $q=p, p+1(=2p-1-(p-2))$, there is an isomorphism
\begin{align} \label{eq:cohrelhcp-2}
H^{q}(X,\overline{HC}^{(p-2)}_{p-2}(O_{X_{A}})) \cong H^{q}(X,\overline{HH}^{(p-2)}_{p-2}(O_{X_{A}}))=0.
\end{align}where the second identity is from Lemma \ref{l:cohrelhhzero} again (let $i=p-2$).

Let $n=p-1$ and $q=p$ in the sequence \eqref{l:lexhhhc}, it has the form
\begin{align*} \label{eq:lexhhhc6}
& \cdots \to H^{p}(X,\overline{HC}^{(p-2)}_{p-2}(O_{X_{A}})) \to H^{p}(X,\overline{HH}^{(p-1)}_{p-1}(O_{X_{A}}))  \\
& \to H^{p}(X,\overline{HC}^{(p-1)}_{p-1}(O_{X_{A}})) \to H^{p+1}(X,\overline{HC}^{(p-2)}_{p-2}(O_{X_{A}})) \to  \cdots.  \notag
\end{align*}The following lemma could be obtained by substituting \eqref{eq:cohrelhcp-2} into this sequence.
\begin{lemma} \label{l:cohrelhcp-1}
In Setting \ref{s:set}, for $A \in Art_{k}$ which is also a graded $k$-algebra $A=\oplus_{l \geq 0}A_{l}$ with $A_{0}=k$, there is an isomorphism
\begin{align*}
H^{p}(X, \overline{HC}^{(p-1)}_{p-1}(O_{X_{A}}))\cong H^{p}(X, \overline{HH}^{(p-1)}_{p-1}(O_{X_{A}})).
\end{align*}
\end{lemma}

Combining this lemma with Lemma \ref{l:iso2} and Lemma \ref{l:cohrelhh}, we obtain the main result of this paper.
\begin{theorem} \label{t:yang1}
In Setting \ref{s:set}, for $A \in Art_{k}$ which is also a graded $k$-algebra $A=\oplus_{l \geq 0}A_{l}$ with $A_{0}=k$, there is an isomorphism
\[
\widehat{CH}^{p}(A) \cong H^{p}(X, \Omega^{p-1}_{X/ k})\otimes_{k}m_{A}.
\]

\end{theorem}

\textbf{Acknowledgments}. The author is very grateful to Spencer Bloch \cite{Bl5} and Jan Stienstra \cite{Stien2} for sharing their ideas. He thanks Phillip Griffiths, Jerome William Hoffman, Kefeng Liu, Yang Shen, Wei Xia, Chenglong Yu and Zhiwei Zheng for discussions and /or comments.

This work is partially supported by the Key Program and the Major Program of Natural Science Research Foundation of Anhui Provincial Education Department (No.2024AH051403 and No.2023AH040225), the Research and Innovation Team of Anhui Province (No.2024AH010023), and the Start-up Research Fund of Chuzhou University (No.2024qd53).

\end{document}